\documentclass[12pt]{amsart}
\usepackage{amsfonts}
\usepackage{amssymb}
\theoremstyle{definition}

\theoremstyle{remark}

\newcommand{\ds}{\displaystyle}

\begin{document}

\centerline{\bf Comptes rendus de l'Academie bulgare des Sciences}

\centerline{\it Tome 35, No 7, 1982}

\vspace{0.4in}

\centerline{\large\bf SUR LE TH\'EOR\`EME DE F. SCHUR POUR UNE VARI\'ET\'E }

\vspace{0.1in}
\centerline{\large\bf PRESQUE HERMITIENNE}

\vspace{0.6in}
\begin{flushright}
{\it MATH\'EMATIQUES
\\ G\'eometrie diff\'erentielle}
\end{flushright}

\vspace{0.3in}
\centerline{\bf  O. T. Kassabov}

\vspace{0.2in}
\centerline{\it (Pr\'esent\'ee par B. Petkanchin, membre de l'Acad\'emie, le 23 mars 1982)}

\vspace{0.2in}
{\bf 1. Introduction.} Soit $M$ une vari\'et\'e riemannienne de tenseur m\'etrique $g$
et tenseur de la courbure $R$. On definit la courbure sectionnelle $K(\alpha)$ d'un
2-plan $\alpha$ par
$$
	K(\alpha)=R(x,y,y,x) \ ,
$$
o\`u $\{ x,y \}$ est une base orthonorm\'ee de $\alpha$. Le th\'eor\`eme de Schur suivant
est bien connu:

{\it Soit $M$ une vari\'et\'e riemannienne de dimension $m>2$, telle que pour chaque point
$p\in M$ la courbure d'un 2-plan arbitraire $\alpha$ dans $T_p(M)$ ne d\'epend pas de
$\alpha$: $K(\alpha)=c(p)$. Alors $c(p)$ ne d\'epend pas de $p$, c'est-\`a dire $M$ est 
\`a courbure sectionnelle constante. }

Soit $M$ une vari\'et\'e presque hermitienne de tenseur m\'etrique $g$, structure presque
complexe $J$ et tenseur de la courbure $R$. Soit $\alpha$ un 2-plan dans $T_p(M)$. 
L'angle $\ds\sphericalangle (\alpha,J\alpha) \in [0,\frac{\pi}2]$ entre $\alpha$ et $J\alpha$ est
donn\'e par
$$
	\cos\sphericalangle(\alpha,J\alpha)=|g(x,Jy)| \ ,
$$
o\`u $\{ x,y \}$ est une base orthonorm\'ee de $\alpha$. Un 2-plan $\alpha$ est dit
$\theta$-holomorphe, si $\sphericalangle(\alpha,J\alpha)=\theta$. Un 0-holomorphe 
(resp.\,$\ds\frac{\pi}2$)-holomorphe) 2-plan est dit aussi holomorphe (resp. anti\-holomorphe). 

On dit que $M$ est \`a courbure sectionnelle $\theta$-holomorphe constante par points, si 
pour chaque point $p\in M$ la courbure d'un 2-plan arbitraire $\theta$-holomorphe 
$\alpha \subset T_p(M)$ ne d\'epend pas de $\alpha$: $K(\alpha)=c(p)$. Particuli\`erement,
si $c(p)$ ne d\'epend pas de $p$, $M$ est dite vari\'et\'e \`a courbure sectionnelle 
$\theta$-holomorphe constante.

Une vari\'et\'e presque hermitienne $M$ est dite $RK$-vari\'et\'e, si
$$
	R(x,y,z,u)=R(Jx,Jy,Jz,Ju)
$$
pour tous $x,\,y,\,z,\,u \in T_p(M)$, $p\in M$. Soit $\nabla$ la connection riemannienne
de $M$. Si $\nabla J=0$, $M$ est dite vari\'et\'e k\"ahlerienne.

En liaison avec le th\'eor\`eme de Schur on peux faire la conjecture suivante:

{\it Soit $\theta \in \ds [0,\frac{\pi}2]$ et soit $M$ une vari\'et\'e presque hermitienne
\`a courbure sectionnelle $\theta$-holomorphe constante par points. Alors $M$ est \`a
courbure sectionnelle $\theta$-holomorphe constante.}

Dans \cite{GV} on a prouv\'e, que c'est vrai pour certains vari\'et\'es presque hermitiennes
en cas o\`u $\theta=0$. Dans section 2 nous allons examiner le cas $\ds\theta\in(0,\frac{\pi}2)$
et dans section 3 - le cas $\theta=\ds\frac{\pi}2$ pour une $RK$-vari\'et\'e.

\vspace{0.3cm}
{\bf 2. Le cas $\ds\theta\in (0,\frac{\pi}2)$.}

\vspace{0.1cm}
{\bf Th\'eor\`eme 1.} Soit $M$ une vari\'et\'e presque hermitienne de dimension $2m \ge 6$
et soit $\ds\theta\in (0,\frac{\pi}2)$. Si $M$ est \`a courbure sectionnelle $\theta$-holomorphe
constante par points, $M$ est \`a courbure sectionnelle constante ou bien $M$ est une vari\'et\'e
k\"ahlerienne \`a courbure sectionnelle $\varphi$-holomorphe constante pour chaque
$\ds\varphi\in (0,\frac{\pi}2)$.

{\bf Preuve.} Si $p\in M$ et $x,\,y$ sont des vecteurs uniques dans $T_p(M)$ avec
$g(x,y)=g(x,Jy)=0$, le plan avec la base $\{ Jx,x\cos\theta+y\sin\theta \}$ est
$\theta$-holomorphe et alors
$$
	R(Jx,x\cos\theta+y\sin\theta,x\cos\theta+y\sin\theta,Jx)=c(p)  \leqno (1)
$$
o\`u $c(p)$ ne d\'epend pas de $x,\,y$. On d\'eduit de (1) et $\ds\theta\in (0,\frac{\pi}2)$
$$
	R(x,Jx,Jx,y)=0 \ ,   \leqno (2)
$$
$$
	R(x,Jx,Jx,x)\cos^2\theta+R(Jx,y,y,Jx)\sin^2\theta=c(p) \ .   \leqno (3)
$$
D'apr\`es (2) et \cite{OK} $M$ est une vari\'et\'e \`a courbure sectionnelle holomorphe 
constante par point $\mu$:
$$
	R(x,Jx,Jx,x)=\mu(p)  \ .   \leqno(4)   
$$
Il r\'esulte de (3) et (4) que $M$ est \`a courbure sectionnelle antiholomorphe constante
par points $\nu$. Donc le tenseur de la courbure a la forme \cite{GG}
$$
	R=\nu\pi_1+\frac{\mu-\nu}3 \pi_2     \leqno (5)
$$
avec 
$$
	\pi_1(x,y,z,u)=g(x,u)g(y,z)-g(x,z)g(y,u) \ ,
$$
$$
	\pi_2(x,y,z,u)=g(x,Ju)g(y,Jz)-g(x,Jz)g(y,Ju)-2g(x,Jy)g(z,Ju) \ .
$$
D'apr\`es \cite{TrV} si $R$ a la forme (5), $M$ est \`a courbure sectionnelle constante $\nu$
ou bien $M$ est une vari\'et\'e k\"ahlerienne \`a courbure sectionnelle holomorphe constante.

\vspace{0.3cm}
{\bf 3. $RK$-vari\'et\'es \`a courbure sectionnelle antiholomorphe constante par points.}

\vspace{0.1cm}
{\bf Th\'eor\`eme 2.} Soit $M$ une $RK$-vari\'et\'e de dimension $2m \ge 6$ \`a courbure
sectionnelle antiholomorphe constante par points. Alors $M$ est \`a courbure sectionnelle 
antiholomorphe constante.

{\bf Preuve.} Nous avons prouv\'e dans \cite{GK}, que si $M$ est une $RK$-vari\'et\'e \`a
courbure sectionnelle antiholomorphe constante par poits $\nu$, alors
$$
	R=\frac16\psi(S) +\nu\pi_1 -\frac{2m-1}3\nu\pi_2  \ ,   \leqno (6)
$$
avec
$$
	\begin{array}{r}\vspace{0.2cm}
		\psi(S)(x,y,z,u)=g(x,Ju)S(y,Jz)-g(x,Jz)S(y,Ju)-2g(x,Jy)S(z,Ju)  \ \ \\
		                 +g(y,Jz)S(x,Ju)-g(y,Ju)S(x,Jz)-2g(z,Ju)S(x,Jy) \ .
	\end{array}
$$
Ici $S$ est le tenseur de Ricci pour $M$.

Pour un point $p\in M$ soit $\alpha$ un 3-plan antiholomorphe dans $T_p(M)$ avec
une base orthonorm\'ee $\{ x,y,z \}$. D'apr\`es la seconde identit\'e de Bianchi
$$
	(\nabla_xR)(y,z,z,y)+(\nabla_yR)(z,x,z,y)+(\nabla_zR)(x,y,z,y)=0  \ .  \leqno (7)
$$ 
On d\'eduit de (6) et (7)
$$
	\begin{array}{r}\vspace{0.2cm}
		2g(y,(\nabla_xJ)z)S(y,Jz)+g(z,(\nabla_yJ)y)S(x,Jz)+g(x,(\nabla_yJ)z)S(z,Jy) \ \ \\  
		          +g(x,(\nabla_zJ)y)S(y,Jz)+g(y,(\nabla_zJ)z)S(x,Jy)+3x(\nu)=0 \ .
	\end{array}    \leqno (8)
$$
Soit $\{e_i,\,Je_i;\,i=1,\hdots,m\}$ une base orthonorm\'ee de $T_p(M)$, telle que 
$S(e_i)=\lambda_ie_i$, $i=1,\hdots,m$. Il suite de (8) pour $x=e_i$, $y=e_j$,
$z=e_k$ $(i\ne j\ne k\ne i)$
$$
	e_i(\nu)=0 \ .  \leqno (9) 
$$ 
On d\'eduit de (9) que $\nu$ ne d\'epend pas du point $p$.

\vspace {0.6cm}
\begin{flushright}
{\it Institute de  Math\'ematiques  \\
Acad\'emie bulgare des Sciences \\
Sofia, Bulgaria}
\end{flushright}

\end{document}